\documentclass [11pt]{article}
\usepackage{amsmath, amsthm, amssymb}
\usepackage{graphicx}
\newfont{\bbb} {msbm10}
\newcommand{\R}{\Bbb{R}}
\newcommand{\bS}{\Bbb{S}}
\newcommand{\bB}{\Bbb{B}}
\newcommand{\T}{\Bbb{T}}

\newcommand{\sbs}{\subset}
\newcommand{\ra}{\rightarrow}

\newcommand{\bg}{{\bar{g}}}

\newcommand{\barB}{{\bar{B}}}

\newcommand{\p}{\partial}

\newcommand{\cW}{{\cal{W}}}

\newcommand{\B}{\Bbb{B}}

\newcommand{\HH}{\Bbb{H}}

\newcommand{\ua}{\underline{a}}
\newcommand{\ub}{\underline{b}}

\newcommand{\0}[1]{_{_{#1}}}

\usepackage[all]{xy}

\begin{document}

\title{Deforming an $\epsilon$-Close to Hyperbolic Metric to a Warped Product}
\author{Pedro Ontaneda\thanks{The author was
partially supported by a NSF grant.}}
\date{}

\maketitle

\begin{abstract} We show how to deform a metric of the form
$g=g_r+dr^2$ to a warped product $\cW g=\sinh^2(r) \,g'+dr^2$ ($g'$ does not depend on $r$),
for $r$ less than some fixed $r\0{0}$.
Our main result establishes to what
extent the {\it warp forced metric} $\cW g$ is {\it close to being
hyperbolic}, if we assume $g$ to be close to hyperbolic.
\end{abstract}

\noindent {\bf \large  Section 0. Introduction.}

First we introduce some notation. The canonical flat metric on
$\R^k$ and the round metric on $\bS^k$ will be denoted by $\sigma\0{\R^k}$
and $\sigma\0{\bS^k}$ respectively.
Let $(M^n, g)$ be a complete Riemannian manifold {\it with center $o\in M$}, that is, the exponential map $\exp_o:T_oM\ra M$ is a diffeomorphism.
Using the exponential map $\exp_o$ we shall sometimes identify $M$ with
$\R^n$, thus we
can write the metric $g$ on $M-\{ o\}=\bS^{n-1}\times\R^+$ as \,\,$g=g_r+dr^2$, where $r$
is the distance to $o$.  The open ball of radius $r$ in $M$, centered at $o$, will be denoted by $B_r=B_r(M)$, and the closed ball by $\barB_r$.
We fix a function $\rho:\R\ra[0,1]$ 
with $\rho (t)=0$ for $t\leq 0$, $\rho(t)=1$ for $t\geq 1$, and $\rho$ constant near 0 and 1.  
\vspace{.1in}

Let $M$ have center $o$ and metric $g=g\0{r}+dr^2$. Fix $r\0{0}>0$.
We define the metric $\bg\0{r\0{0}}$ on $M-\{o\}$ by:
\begin{center}$\bg \0{r\0{0}}\,=\, \sinh^2\, (r)
\big(  \frac{1}{\sinh^2(r\0{0})}\big) g\0{r\0{0}}\, +\, dr^2$. \end{center}
Note that this metric is a warped product (warped by $\sinh$). 
Note also that to define $\bg\0{r\0{0}}$ we are using the identification
$M-\{ o\}=\bS^{n-1}\times\R^+$ given by the original metric $g$. 
We now force the metric $g$ to be equal to $\bg\0{r\0{0}}$ on  $\barB_{r\0{0}}=\barB_{r\0{0}}(M)$ and stay equal to $g$ outside $B_{r\0{0} +\frac{1}{2}}$.
For this we define the {\it warp forced  (on $B_{r\0{0}}$) metric } as:

\begin{center}$
\cW\0{r\0{0}}\, g\, =\, \rho\0{r\0{0}}\, \bg\0{r\0{0}}\, +\, 
(1-\rho\0{r\0{0}})\, g.
$\end{center}

\noindent where $\rho\0{r\0{0}}(t)=\rho(2t-2r\0{0})$.
Hence we have\\

\noindent{\bf (0.1)}\hspace{1.2in}
 $\cW\0{r\0{0}} g\,=\,\left\{ \begin{array}{lllll}
\bg\0{r\0{0}}&& {\mbox{on}}& & \barB_{r\0{0}}\\\\
g&&{\mbox{outside}}& & B_{r\0{0}+\frac{1}{2}}.
\end{array}\right.$\\

We call the process $g\mapsto\cW\0{r\0{0}} g$ {\it warp forcing}.
Note that if we choose $g$ to be the warped-by-sinh hyperbolic metric $g = \sinh^ 2(t) \sigma\0{\bS^{n-1}} + dt^2$, then $\cW\0{r\0{0}} g = g$. This suggests that
if $g$ is in some sense close to being hyperbolic, then
$\cW\0{r\0{0}} g$ should also be close to hyperbolic.
The purpose of this paper is to quantify this last
 statement, that is, to answer the following question:  if $g$ is {\it $\epsilon$-close to a hyperbolic metric} then to what extent is
the warp forced metric  $\cW\0{r\0{0}} g$ close to hyperbolic?
The answer is that $\cW\0{r\0{0}} g$ is $\eta$-close to hyperbolic
where $\eta$ depends on $\epsilon$ and $r\0{0}$.
The term ``$\epsilon$-close to a hyperbolic metric"
used above refers to a chart-by-chart concept; it is introduced in the next paragraph.
\vspace{.1in}

Let $\bB$ be the unit open $(n-1)$-ball 
with the flat metric $\sigma\0{\R^{n-1}}$.  Write  $I_\xi=(-1-\xi,1+\xi)$, $\xi\geq 0$.
Our {\it basic models} are $\T_\xi=\bB\times I_\xi$,  with hyperbolic metric $\sigma=e^{2t}\sigma\0{\R^{n-1}}+dt^2$. 
The number $\xi$ is called the {\it excess} of $\T_\xi$.
(The reason for introducing $\xi$ will become clear in 
the Main Theorem below; see also the remark after the Theorem).
Let $(M,g)$ be a Riemannian manifold and $S\sbs M$. We say that $g$ is $\epsilon$-{\it close to hyperbolic on $S$} if there is $\xi\geq 0$ such that for every $p\in S$ there is an {\it $\epsilon$-close to hyperbolic 
chart with center $p$}, that is, there is a chart
$\phi :\T_\xi\ra M$, $\phi(0,0)=p$,  such that 
$|\phi^*g-\sigma|\0{C^2}<\epsilon$.  
The number $\xi$ is called the {\it excess } of the charts.
We stress that $\xi$ is independent of $p$. 
Here $|.|\0{C^2}$ is the $C^2$-norm
(see Section 1).
\vspace{.1in}

Let $(M,g)$ have center $o$ and $S\sbs M$. We say that $g$ is {\it radially $\epsilon$-close to hyperbolic on $S$ (with respect to $o$)} if, 
for every $p\in S$ there is an  $\epsilon$-close to hyperbolic 
chart $\phi$ with center $p$ and, in addition, the chart 
$\phi$ respects the product structure of $\T_\xi$
and $M-o=\bS^{n-1}\times\R^+$, that is $\phi(. , t)=(\phi\0{1}(.), t+a)$, 
where the constant 
 $a$ depends on 
$\phi$, and $\phi\0{1}$ is some function independent of $t$ (equivalently, $\phi_1$ is a chart on $M$). Here the ``radial" directions
are $(-1-\xi,1+\xi)$ and $\R^+$ in $\T_\xi$ and $M-o$, respectively.
\vspace{.1in}

As mentioned before, our main result below shows that if $g$ is radially $\epsilon$-close to hyperbolic then
the warp forced metric  $\cW\0{r\0{0}} g$ is radially $\eta$-close to hyperbolic,
where $\eta$ depends on $\epsilon$ and $r\0{0}$. 
In the next Theorem we assume $\xi>1$
and $r\0{0}\geq 3+2\xi$. \vspace{.1in}

\noindent {\bf  Theorem.} {\it  Let $(M,g)$ have center $o$, and $S\sbs M$. If $g$ is radially $\epsilon$-close to hyperbolic on $S$, with charts of excess $\xi$,
then $\cW\0{r\0{0}} g$ is radially $\eta$-close to hyperbolic on $S-\barB\0{r\0{0}-1-\xi}$ with charts of excess $\xi-1$, provided $\eta\geq e^{27+12\xi}\big(e^{-2r\0{0}}+\epsilon\big)$.}
\vspace{.1in}

\noindent{\bf Remark.} Note that warp forcing reduces the excess of the charts by 1.
This was one of the motivations to introduce the excess $\xi$.
\vspace{.1in}

The results in this paper are used to construct negatively curved Riemannian smoothings of Charney-Davis strict
hyperbolizations of manifolds \cite{ChD}, \cite{O}. In the next paragraph we give an idea how the Theorem in this paper is used in \cite{O}.
\vspace{.1in}

In the same way that a cubical complex is made of basic pieces (the cubes
$\square^k$),
the hyperbolization $h(K)$ of a cubical complex $K$ is 
also made of
basic pieces: pre-fixed hyperbolization pieces $X^k$. Indeed one begins with a cubical complex $K$ and replaces each cube of dimension $k$
by the hyperbolization piece of the same dimension. Cube complexes have
a piecewise flat metric induced from the flat geometry of the cubes.
Likewise the Charney-Davis hyperbolizations have a piecewise hyperbolic
structure because the Charney-Davis hyperbolization pieces
are hyperbolic manifolds (compact, with boundary and corners).
To see how singularities appear one can first think about the manifold 2-dimensional 
cube case. If $K^2$ is a 2-dimensional manifold cube complex then
its piecewise flat metric is Riemannian outside the vertices. A vertex is
a singularity if and only if the vertex does not meet exactly four cubes.
The picture is exactly the same for $h(K^2)$. 
These point singularities in $h(K^2)$ can be smoothed out
easily using warping methods.
In higher dimensions
the singularities of $K^n$ and $h(K)$ appear in (possibly the whole of)
the codimension 2 skeletons $K^{(n-2)}$ and $h(K^{(n-2)})$, respectively.
In \cite{O} the idea of smoothing the piecewise hyperbolic metric on
$h(K)$ is to do it inductively down the dimension of the skeleta.
One begins with the $(n-2)$-dimensional pieces $X^{n-2}$. Transversally
to each $X^{n-2}$ (that is, on the union of geodesic segments emanating
perpendicularly to $X^{n-2}$, from a fixed point in $X^{n-2}$) one has
essentially the 2-dimensional picture mentioned above. Once we
solve this transversal problem we extend this transversal smoothing by taking a warp product
with $X^{n-2}$; we called this product method {\it hyperbolic extension}
\cite{O2}.
This gives a smoothing on a (tubular) neighborhood of the piece $X^{n-2}$.
Caveat: we do not want to actually have a smoothing on a neighborhood
of the whole of $X^{n-2}$, since we will certainly have 
matching problems for different $X^{n-2}$ meeting on a common
$X^{n-3}$; so we only want a smoothing on a neighborhood of the $Z^{n-2}$,
where $Z^{n-2}\sbs X^{n-2}$ is just a bit ``smaller" than $X^{n-2}$,
so that the neighborhoods of the $Z^{n-2}$ are all disjoint. Next 
step is to smooth around the $X^{n-3}$ (or, specifically the $Z^{n-3}$).
The metric is already smooth outside a neighborhood of the $(n-3)$-skeleton. Transversally to each $X^{n-3}$ we
have a 3 dimensional problem.
(It helps to have a 3 dimensional picture in mind, like in dimension 2).
It happens that if we did things with care in the first step (around the
$Z^{n-2}$) the metric in the 3 dimensional transversal problem
is radially $\epsilon$-close to hyperbolic outside some large ball B.
If this metric was a warped product we could use the {\it two variable
warping deformation} given in \cite{O2} to extend the metric to a Riemannian metric
on the ball B, getting rid, in this way, of the transverse
singularity. But
the metric in the 3 dimensional transversal problem
is not warped, hence the need for the Theorem in this paper:
one takes a radially $\epsilon$-close to hyperbolic metric
and deforms it to a warped metric inside a ball, and the resulting
metric is still radially $\eta$-close to hyperbolic, with an $\eta$
that can be controlled. Once the transversal 3 dimensional problem
is solved we extend this smoothing to neighborhoods of the $Z^{n-3}$
using hyperbolic extension. Next we do the same for the $Z^{n-4}$
and so on. About the excess: since warp forcing reduces the excess
by 1, one begins with a large excess at codimension 2, so that when
we arrive at codimension $n$ one still has positive excess;
therefore in the Theorem above one
should think of the $\xi$ as fixed, while the $r\0{0}$ as being
as large as wanted, $\epsilon$ as small as desired, and the set $S$ as the complement of the
ball of radius $r\0{0}-1-\xi$.\vspace{.1in}

In Section 1 we give some definitions and a useful lemma. In Section 2
we give some estimates on changing warping functions. In Section 3
we do warp forcing locally. In Section 4 prove the Theorem.
\vspace{.1in}

We are grateful to the referee for the many comments and suggestions.\\

\noindent {\bf \large  Section 1. Preliminaries.}\\ Let $A\sbs \R^n$ be an open set.
Let $|.|\0{C^2(A)}$ denote the uniform $C^2$-norm of $\R^l$-valued functions on $A$, i.e. if $f=(f\0{1},...,f\0{l}):A\ra \R^l$,
then $|f|\0{C^2(A)}=sup\0{z\in A,
\,\,1\leq i\leq l,\,\,1\leq j,k\leq n}\{ |f\0{i}(z)|, |\p\0{j}f\0{i}(z)|, |\p\0{j,k}f\0{i}(z)|\}$. Sometimes we will write $|.|\0{C^2}=|.|\0{C^2(A)}$ when the context is clear.
Given a Riemannian metric $g$ on $A$, the number $|g|\0{C^2(A)}$ is computed considering $g$ as the $\R^{n^2}$-valued function $z\mapsto (g_{ij}(z))$ where, as usual,
$g\0{ij}=g(e_i,e_j)$, and the $e_i$'s are the canonical vectors in $\R^{n}$. 
\vspace{.1in}

The $C^2$-norm $|.|\0{C^2}$ mentioned in the  
definition of an  $\epsilon$-close to hyperbolic Riemannian manifold in the Introduction is $|.|\0{C^2}=|.|\0{C^2(\T_\xi)}$. 
If $(M,g)$ is $\epsilon$-close to hyperbolic (or radially $\epsilon$-close 
to hyperbolic) we will also say that the metric $g$ is $\epsilon$-close to hyperbolic (or radially $\epsilon$-close to hyperbolic).
\vspace{.1in}

Note that for the metric $\sigma=e^{2t}\sigma\0{\R^{n-1}}+dt^2$ on our model 
$\T_\xi$ we have $|\sigma|\0{C^2(\T_\xi)}=4e^{2+2\xi}$.\vspace{.1in}

\noindent {\bf Remarks.} \\
\noindent {\bf 1.} The definition of radially $\epsilon$-close to hyperbolic metrics 
is well-suited to studying metrics of the form $g\0{t}+dt^2$ for $t$ large, but for
small $t$ this definition has some drawbacks because: (1) we need some space to fit the charts, and (2) the form of our specific fixed model $\T_\xi$.
An undesired consequence is that punctured hyperbolic space $\HH^n-\{ o\}=\bS^{n-1}\times\R^+$
(with warped product $\sinh^2(t)\sigma\0{\bS^{n-1}}+dt^2$) is not radially $\epsilon$-close to hyperbolic
for $t$ small.

\noindent{\bf 2.} 
In \cite{O} we actually need warped metrics with warping functions that are multiples of hyperbolic
functions. All these functions are close to the exponential $e^t$ (for $t$ large), so instead of introducing one model for each hyperbolic function
we introduced only the exponential model. In the next section we show
the effect of changing warping functions.\vspace{.1in}

We will need the following lemma.\vspace{.1in}

\noindent {\bf Lemma 1.1.} {\it Let $g\0{i}$ be metrics on $\T_\xi $ such that
$|g\0{i}-\sigma|\0{C^2(\T_\xi)}<\epsilon\0{i}$, for $i=1,\,2$. Let $\lambda:\T_\xi\ra [0,1]$ be smooth with $|\lambda|\0{C^2(\T_\xi)}$
finite. Then}
 $$\Big|\,\,\,\,\lambda\, g\0{1}\,+\, (1-\lambda)\,g\0{2}\,\,\,\,\,-\,\,\,\,\,\sigma \,\,\,\,\Big|_{C^2(\T_\xi)}
 \,\,\,\,\,<\,\,\,\,\,4\,\,\,(1+|\lambda|\0{C^2(\T_\xi)})\,\,\,
 (\epsilon\0{1}+\epsilon\0{2}).$$

\noindent {\bf Proof.} The proof follows from the triangle inequality,
Leibniz rule and the equality $\big(\lambda\, g\0{1}\,+\, (1-\lambda)\,g\0{2}\big)-\sigma=
\lambda\,(g\0{1}-\sigma)+(1-\lambda)(g\0{2}-\sigma)$.
This proves the lemma.\vspace{.2in}

\noindent {\bf \large Section 2.  Warping with  $\sinh\, t$.}\\
The metric of our basic hyperbolic model $\T_\xi$ is an exponentially warped metric. Here we show that we can 
change the exponential by multiples of $\sinh(t)$, for $t$ large.\vspace{.1in}

In what follows we will often consider metrics $h$ on
$\T_\xi$ of the form $h=h\0{t}+dt^2$.
Recall $I_\xi=(-1-\xi,\, 1+\xi)$. 
\vspace{.1in}

\noindent {\bf Lemma 2.1.} 
{\it  For ${\sf{r}}\geq 2+\xi$ \, we have $\big| e^{-2t}\big(\frac{\sinh\, (t+{\sf{r}})}{\sinh\, ({\sf{r}})}\big)^2\, -\, 1   \big|\0{C^2(I_\xi)}\, <\,43\,e^{2+2\xi}\,e^{-2{\sf{r}}}.$ }\vspace{.1in}

\noindent {\bf Proof.} Write 
$e^{-t}\,\frac{\sinh\, (t+{\sf{r}})}{\sinh\, ({\sf{r}})}
\, -\, 1 \,=\,\frac{1-e^{-2t}}{1-e^{-2{\sf{r}}}}
\,e^{-2{\sf{r}}}$.
Since ${\sf{r}}\geq 2$, we have $\frac{1}{1-e^{-2{\sf{r}}}}
\leq \frac{1}{1-e^{-4}}<1.02$.
Differentiating $\big(\frac{1}{1-e^{-2{\sf{r}}}}\big)(1-e^{-2t})\,e^{-2{\sf{r}}}$ twice,
together with the previous two facts give
the following estimate:

{\small
$$\,\,\big|e^{-t}\bigg(\frac{\sinh\, (t+{\sf{r}})}{\sinh\, ({\sf{r}})}\bigg)\, -\, 1   \big|_{C^2(I_\xi)}\,\,\, < \, \,\,(1.02)\,\, (4e^{2+2\xi})\,\,e^{-2{\sf{r}}}\,\,\,=\,\,\,4.08 \,\,e^{2+2\xi} 
\,\,e^{-2{\sf{r}}}.$$}

This estimate together with the triangle inequality and the
hypothesis ${\sf{r}}\geq 2+\xi$    give the
following estimate:

{\small $$\,\,\big| e^{-t}\bigg(\frac{\sinh\, (t+{\sf{r}})}{\sinh\, ({\sf{r}})}\bigg)\,+1\, \big|_{C^2(I_\xi)}\, \leq 2\,+\,4.08\,e^{2+2\xi}e^{-2{\sf{r}}}\,=\,  2\,+\,4.08\,e^{2+2\xi-2{\sf{r}}}\,\leq\,
 \,2+\,4.08\,e^{-2}\,<\, 2.6.\,$$}

To prove the lemma write
\begin{center}
{\small $e^{-2t}\big(\frac{\sinh\, (t+{\sf{r}})}{\sinh\, ({\sf{r}})}\big)^2\, -\, 1 =
\bigg(e^{-t}\big(\frac{\sinh\, (t+{\sf{r}})}{\sinh\, ({\sf{r}})}\big)
\, -\, 1\bigg)\bigg( e^{-t}\big(\frac{\sinh\, (t+{\sf{r}})}{\sinh\, ({\sf{r}})}
\big)\, +\, 1\bigg). $}\end{center} 
\noindent This together with the previous two estimates
and the Leibniz rule give

{\small $$\big|\,e^{-2t}\big(\frac{\sinh\, (t+{\sf{r}})}{\sinh\, ({\sf{r}})}\big)^2\, -\, 1\,\big|
_{C^2(I_\xi)} \,\,\leq\,\,4\,(4.08e^{2+2\xi}e^{-2{\sf{r}}})\,2.6\,\,<\,\,
43\,e^{2+2\xi}\,e^{-2{\sf{r}}}.$$}
This proves the lemma. \vspace{.2in}

 Let $\nu:I_\xi\ra \R^+$ be smooth. For a metric $f=f\0{t} + dt^2$ on $\T_\xi$ we write
$f^{\nu}=\nu f\0{t}+dt^2$. \vspace{.1in}

\noindent {\bf Lemma 2.2.} {\it We have $\big| f^{\nu} -f \big|_{C^2(\T_\xi)}\, \leq\, 4\, \big| \nu-1   \big|_{C^2(I_\xi)}\, \big| f\big|_{C^2(\T_\xi)}$.}\vspace{.1in}

\noindent {\bf Proof.} Just note that $f^{\nu}-f=(\nu-1)f\0{t}$ and differentiate twice. This proves the lemma. \vspace{.1in}

Recall that the metric on our model $\T_\xi$ is $\sigma=e^{2t}\sigma\0{\R^{n-1}}+dt^2$. 
\vspace{.1in}

\noindent {\bf Lemma 2.3.} {\it Let $f=f\0{t}+dt^2$ be
a metric on $\T_\xi$ such that
 $|f-\sigma|\0{C^2(\T_\xi)}<\epsilon$. 
Let $\nu=e^{-2t}\big(\frac{\sinh\, (t+{\sf{r}})}{\sinh\, ({\sf{r}})}\big)^2$. 
Assume ${\sf{r}}\geq 2+\xi$. Then }
\begin{enumerate}
\item[{\it (1)}] $ \big | f^\nu-f \big |\0{C^2(\T_\xi)}\, <\,172\, e^{2+2\xi}\, (\epsilon +4e^{2+2\xi})\, e^{-2{\sf{r}}}$.
\item[{\it (2)}] $|f^\nu-\sigma|\0{C^2(\T_\xi)}<  \, 688\,  e^{4+4\xi}\, \big(\epsilon+ e^{-2{\sf{r}}}\,\big)$.
\end{enumerate}
\noindent {\bf Proof.} Item 1
follows from 2.1,  2.2, and the fact that $|f|\0{C^2(\T_\xi)}\leq |f-\sigma|\0{C^2(\T_\xi)}+|\sigma|\0{C^2(\T_\xi)}<\epsilon+4\,e^{2+2\xi}$. To prove 
item 2 note that it follows from item 1 and the hypothesis
$|f-\sigma|\0{C^2(\T_\xi)}<\epsilon$ that\\

\hspace{.2in}$|f^\nu-\sigma|\0{C^2(\T_\xi)}\leq |f-\sigma|\0{C^2(\T_\xi)}+|f^\nu-f|\0{C^2(\T_\xi)}
<\,\epsilon\,+\,\,172\, e^{2+2\xi}\, (\epsilon +4e^{2+2\xi})\, e^{-2{\sf{r}}}$\\

\hspace{2.98in}$=\,\,(1+172\,e^{2+2\xi-2{\sf{r}}})\,\epsilon\,\,+\,\,172\,e^{2+2\xi}\,4\,e^{2+2\xi}\, e^{-2{\sf{r}}}$\\

\hspace{2.98in}$<\,\,172\,e^{2+2\xi}\,4\,e^{2+2\xi}\,(\epsilon+ e^{-2{\sf{r}}})$\\

\hspace{2.98in}$=\,\,688\,e^{4+4\xi}\,(\epsilon+ e^{-2{\sf{r}}}).$\vspace{.1in}

This proves the lemma. \vspace{.1in}

As in Lemma 2.3 let $\nu=e^{-2t}\big(\frac{\sinh\, (t+{\sf{r}})}{\sinh\, ({\sf{r}})}\big)^2$. Lemma 2.1 says that  
$\big| \nu -\, 1   \big|_{C^2(I_\xi)}\, <\,43\,e^{2+2\xi}\,e^{-2{\sf{r}}}.$
\vspace{.1in}

Let $s\in I_\xi$.  Write $\nu\0{s}(t)=\nu(t-s)$
with $\nu$ as above. \vspace{.1in}

\noindent {\bf Lemma 2.4.} 
{\it  For ${\sf{r}}\geq 2+\xi$ and $s\in I_\xi$\, we have $\big| \nu\0{s}\, -\, 1   \big|\0{C^2(I_\xi)}\, <\,43\,e^{4+4\xi}\,e^{-2{\sf{r}}}.$ }\vspace{.1in}

\noindent {\bf Proof.}
For $t\in I_\xi$ we have $t-s\in I_{1+2\xi}$. This together with Lemma 2.1
imply $| \nu\0{s}\, -\, 1 |\0{C^2(I_\xi)}\, <\,43\,e^{2+2(2\xi+1)}\,e^{-2{\sf{r}}}$. 
This proves the lemma.\vspace{.1in}

The next lemma is similar to 2.3, with $\nu\0{s}$ replacing $\nu$ in the
conclusion.\vspace{.1in}

\noindent {\bf Lemma 2.5.} {\it Let $f=f\0{t}+dt^2$ be
a metric on $\T_\xi$ such that
 $|f-\sigma|\0{C^2(\T_\xi)}<\epsilon$. 
Let $\nu=e^{-2t}\big(\frac{\sinh\, (t+{\sf{r}})}{\sinh\, ({\sf{r}})}\big)^2$. 
Assume ${\sf{r}}\geq 2+\xi$. Then }
\begin{enumerate}
\item[{\it (1)}] $ \big | f^{\nu\0{s}}-f \big |\0{C^2(\T_\xi)}\, <\,172\, e^{4+4\xi}\, (\epsilon +4e^{2+2\xi})\, e^{-2{\sf{r}}}$.
\item[{\it (2)}] $|f^{\nu\0{s}}-\sigma|\0{C^2(\T_\xi)}<  \, 688\,  e^{6+6\xi}\, \big(\epsilon+ e^{-2{\sf{r}}}\,\big)$.
\end{enumerate}

The proof is the same as the proof of Lemma 2.3, but uses Lemma 2.4 instead 
of Lemma 2.1.
\vspace{.2in}

\noindent {\bf \large Section 3. Local warp forcing.}  \\
Here we give a kind of a local version to warp forcing. \vspace{.1in}

Let $a$ be a metric on $\B^{n-1}$. For a fixed $s\in I_\xi$ we denote by $\ua_s$ the warped metric $e^{2(t-s)}a+dt^2$ on $\T_\xi=\B^{n-1}\times I_\xi$. \vspace{.1in}

\noindent {\bf Lemma 3.1.} {\it Let $s\in I_\xi$ and let $a$, $b$ be metrics on $\B^{n-1}$ with \,\,$|\, a\, -\, b\, |\0{C^2(\B^{n-1})}\, < \,\epsilon$. Then 
\,$|\, \ua_s\, -\, \ub_s\, |\0{C^2(\T_\xi)}\, <\,16 \, e^{4+4\xi}\,\epsilon$.}

\noindent {\bf Proof.}
Just compute the derivatives of $\ua_s-\ub_s=e^{2(t-s)}(a-b)$. This proves the lemma. 
\vspace{.1in}

The next lemma gives local estimates
(that is, on the model $\T_\xi$) needed for 
global warp forcing
estimates.\vspace{.1in}

\noindent {\bf Lemma 3.2.} {\it Let $h=h\0{t}+dt^2$ be a  metric on $\T_\xi$ with $|h-\sigma|\0{C^2(\T_\xi)}<\epsilon$. 
Fix $s\in I_\xi$ and
consider the warped-by-exponential metric ${\underline{h\0{s}}}_s=e^{2(t-s)}h\0{s}+dt^2$ on $\T_\xi$. Then   $|\underline{h\0{s}}_s-\sigma|\0{C^2(\T_\xi)}<16\,e^{4+4\xi}\,\epsilon$.} \vspace{.1in}
  
\noindent {\bf Proof.} By hypothesis
we have $|\, (h\0{t}+dt^2)\,-\, (e^{2t}\sigma\0{\R^{n-1}}+dt^2)\, |\0{C^2(\T_\xi)}\,<\,\epsilon$. Therefore, taking $t=s$ we get
$|\, h\0{s}\, -\, e^{2s}\sigma\0{\R^{n-1}}\, |\0{C^2(\B^{n-1})}\, < \,\epsilon$. Note that ${\underline{e^{2s}\sigma\0{\R^{n-1}}}}_s=e^{2t}\sigma\0{\R^{n-1}}+dt^2=\sigma$. This together with
Lemma 3.1 implies that $|{\underline{h\0{s}}}_s-\sigma |\0{C^2(\T_\xi)}<16\,e^{4+4\xi}\epsilon$. This completes the proof of Lemma 3.2. 
\vspace{.2in}

\noindent {\bf \large  Section 4.  Proof of the Theorem.}\\
Let $(M^n, g)$ be a complete Riemannian manifold with center $o\in M$. Recall that we
can write the metric on $M-\{ o\}=\bS^{n-1}\times\R^+$ as \,\,$g=g_r+dr^2$. 
Also $B\0{r}$ is the closed ball on $M$ of radius $r$ centered at the center $o$. Let $r_0\geq 3+2\xi$.
We assume that $g$ is radially $\epsilon$-close to hyperbolic on some $S\sbs
M$, with charts of excess $\xi$.
We have to prove that $\cW\0{r\0{0}} g$ is radially $\eta$-close to 
hyperbolic on $S-B\0{r\0{0}-1-\xi}$, with charts of excess $\xi-1$, where 
$\eta=e^{27+12\xi}\big(e^{-2r\0{0}}+\epsilon\big)$.
\vspace{.1in}

Assume $p=(x,{\sf{r}})\in S \sbs\bS^{n-1}\times \R_+=M-\{o\}$ and 
$p\notin \barB\0{r\0{0}-(1+\xi)}$ (equivalently ${\sf{r}}> r\0{0}-(1+\xi)$).  Since the metric $g$ is radially $\epsilon$-close to hyperbolic on $S$, with charts of excess $\xi$, there is a
 radially $\epsilon$-close to hyperbolic chart $\phi:\T_\xi\ra M$
 centered at $p$.
 This means that $\phi(0,0)=p$, $\phi$ is {\it radial}, and
 $|\phi^*g-\sigma|\0{C^2(\T_\xi)}<\epsilon$. 
 Here by {\it radial} we mean that $\phi$ respects product structures (see 
 the definition of a radially $\epsilon$-close to hyperbolic chart in 
 the Introduction).
 To prove Theorem we will prove that
 the restriction $\phi|\0{\T_{\xi-1}}:\T_{\xi-1}\ra M$
 is a radially $\eta$-close to hyperbolic chart for
 $\cW\0{r\0{0}}$ centered at $p$. That is, we will show that
 $|\phi^*\big( \cW\0{r\0{0}}\, g \big )-\sigma|\0{C^2(\T_{\xi-1})}< \eta$.
 We have  three cases. \vspace{.1in}

\noindent {\bf First  case.}  $p\notin B_{r\0{0}+\frac{1}{2}+(1+\xi)}$\\
Then the image of $\phi$ lies outside $ B_{r\0{0}+\frac{1}{2}}$. By (0.1) we have that $\cW\0{r\0{0}}=g$ 
outside $ B_{r\0{0}+\frac{1}{2}}$. Hence the chart $\phi$ is also
a radially $\epsilon$-close to hyperbolic chart for $\cW\0{r\0{0}} g$
centered at $p$ with excess $\xi$.
This shows  the metric 
$\cW\0{r\0{0}} g$ is radially $\epsilon$-close to hyperbolic outside $B_{r\0{0}+\frac{1}{2}+(1+\xi)}$, with charts of excess $\xi$. \vspace{.1in}

\noindent {\bf Second case.}  $p\in B_{r\0{0}+\frac{1}{2}+(1+\xi)}-
B_{r\0{0}+\frac{1}{2}+\xi}$\\
Then the image of the restriction $\phi|\0{\T_{\xi-1}}$
of $\phi$ to $\T_{\xi-1}$ does not intersect $B_{r\0{0}+\frac{1}{2}}$.
Hence as in the first case, by (0.1), the chart $\phi|\0{\T\0{\xi-1}}$ is an
$\epsilon$-close to hyperbolic chart for $\cW\0{r\0{0}} g$,
centered at $p$, but with excess
$\xi-1$. Clearly $\phi|\0{\T\0{\xi-1}}$ is also radial. This shows  the metric 
$\cW\0{r\0{0}} g$ is radially $\epsilon$-close to hyperbolic on $B_{r\0{0}+\frac{1}{2}+(1+\xi)}-B_{r\0{0}+\frac{1}{2}+\xi}$, with charts of excess $\xi-1$.
\vspace{.1in}

\noindent  {\bf Third case.}  $p\in B_{r\0{0}+\frac{1}{2}+\xi}$\\
The condition $p\in B_{r\0{0}+\frac{1}{2}+\xi}$ is equivalent to
${\sf{r}}< r\0{0}+\frac{1}{2}+\xi$. Since by hypothesis $p\notin B\0{r\0{0}-1-\xi}$
we get $r\0{0}-(1+\xi)<{\sf{r}}< r\0{0}+\frac{1}{2}+\xi$. 
Recall $\phi:\T_\xi\ra M $ is a radially $\epsilon$-close to hyperbolic chart of $g$ centered at $p=(x,{\sf{r}})$.
Write $h=\phi^*g$.  Since $\phi$ is radial we have 
that $h$ has the form $h=h\0{t}+dt^2$ with $h\0{t}=\phi^*g\0{t+{\sf{r}}}$.  Moreover

\begin{equation*}
|h-\sigma|\0{C^2(\T_\xi)}\,\,<\,\,\epsilon.
\tag{1}
\end{equation*}

Write $ s=r\0{0}-{\sf{r}}$, thus $-\frac{1}{2}-\xi< s<1+\xi$. 
In particular we have $s\in I_\xi$. Also $h\0{s}=\phi^*g\0{r\0{0}}$.
Recall that in the Introduction we defined the warped product
$\bg \0{r\0{0}}$ as $\bg \0{r\0{0}}=\sinh^2\, (r)\big(  \frac{1}{\sinh^2(r\0{0})}\big) g\0{r\0{0}}\, +\, dr^2$. Since $\phi$ is radial we have 
$\phi^*(\bg \0{r\0{0}})=
\sinh^2\, (t+{\sf{r}})\big(  \frac{1}{\sinh^2(r\0{0})}\big)\,\phi^* g\0{r\0{0}}+dt^2$.
Therefore
\begin{equation*}
\phi^*(\bg \0{r\0{0}})=\,\,\sinh^2\, ((t-s)+r\0{0})
\big(  \frac{1}{\sinh^2(r\0{0})}\big)\,h\0{s}+dt^2.
\tag{2}
\end{equation*}

Note that
$e^{2(t-s)}\nu_s(t)=\frac{\sinh^2 ((t-s)+r\0{0})}{\sinh^2 (r\0{0})}$, 
where $\nu(t)=e^{-2t}\frac{\sinh^2(t+r\0{0})}{\sinh^2(r\0{0})}$ 
and $\nu_s(t)=\nu(t-s)$, as in Section 2.
Using this and the notation in sections 2 and 3, equation (2) 
can be rewritten as

\begin{equation*}
\phi^*\big(\, {\bar{g}}\0{r\0{0}}\,\big)\,\,=\,\,f^{\nu\0{s}}.
\tag{3}
\end{equation*}

\noindent where $f={\underline{h\0{s}}}\0{s}$.  Equation (1) and Lemma 3.2 imply
that $|f-\sigma|\0{C^2(\T_\xi)}<16\,e^{4+4\xi}\epsilon$.
This together with the second item of Lemma 2.5 imply
 
\begin{equation*} |f^{\nu\0{s}} -\sigma|\0{C^2(\T_\xi)}<
688\,e^{6+6\xi}\,\big(e^{-2{\sf{r}}}+16\,e^{4+4\xi}\epsilon\big).\,
\tag{4}
\end{equation*}

\noindent (To apply Lemma 2.5 we need the condition ${\sf{r}}\geq
2+\xi$. This follows from ${\sf{r}}>r\0{0}-(1+\xi)$
and the hypothesis $r\0{0}\geq 3+2\xi$.)
\vspace{.1in}

\noindent From the definition of $\cW\0{r\0{0}}\, g $ given in the Introduction
and the fact that $\phi$ is radial
we have

\begin{equation*}
\phi^*\big( \cW\0{r\0{0}}\, g \big )\,=\,\rho\0{s}\phi^*\big(\bg\0{r\0{0}}\big)
\,+\,\big(1-\rho\0{s}\big) \,\phi^* g\,=\,
\rho\0{s}f^{\nu\0{s}}
\,+\,\big(1-\rho\0{s}\big) \, h.
\tag{5}
\end{equation*}
\vspace{.1in}

\noindent where $\rho\0{s}(t)=\rho(2t-2s)$, and $\rho$ as in the Introduction.
From (5), (4), (1), and Lemma 1.1 we get that 
$|\phi^*\big( \cW\0{r\0{0}}\, g \big )-\sigma|\0{C^2(\T_\xi)}<
\epsilon'$ with

$$
\epsilon'\,\,=\,\,\, 4\, \big( 1+ |\rho\0{s}|\0{C^2(I_\xi)}  \big) 
\Big(\,\,\epsilon\,\,+\,\,688\,e^{6+6\xi}\,\big(e^{-2{\sf{r}}}+16
\,e^{4+4\xi)}\epsilon  \big)\,\,\Big).$$

Note that

$$\epsilon'\,\,<\,\, 4\, \big( 1+ |\rho\0{s}|\0{C^2(I_\xi)}  \big)
\,\, \Big[\,\,1\,\,+\,\,688\,e^{6+6\xi}\,16\,e^{4+4\xi}\,\,
\Big](e^{-2{\sf{r}}}+\epsilon )\,\,<\,\,44033\, 
\big( 1+ |\rho\0{s}|\0{C^2(I_\xi)}  \big)\,e^{10+10\xi}
(e^{-2{\sf{r}}}+\epsilon ).$$

\noindent A calculation shows that we can take $|\rho \0{s}(*)|\0{C^2(I_\xi)}=|\rho (2*)|\0{C^2(\R)}<48$.
This implies that we can take
$\epsilon'<(44033)\, (49)\,e^{10+10\xi}(e^{-2{\sf{r}}}+\epsilon)$.
This together with ${\sf{r}}>r\0{0}-(1+\xi)$ imply that we can take
$\epsilon'<(44033)\, (49)\,e^{12+12\xi}(e^{-2r\0{0}}+\epsilon)= 2157617\,e^{12+12\xi}(e^{-2r\0{0}}+\epsilon)<e^{27+12\xi}(e^{-2r\0{0}}+
\epsilon)$. Note that the excess of the charts in this third case
is also $\xi$. This proves the Theorem.

Pedro Ontaneda

SUNY, Binghamton, N.Y., 13902, U.S.A.

\end{document}